\documentclass[11pt,a4paper]{article}

\usepackage{amsmath} 
\usepackage{amssymb}
\usepackage{dsfont}
\usepackage{graphicx}
\newtheorem{theorem}{Theorem}[section]

\newtheorem{definition}[theorem]{Definition}
\newtheorem{remark}{Remark}

\newcommand{\D}{{\, \rm d}}

\newcommand{\R}{\mathbb{R}}

\newcommand{\bu}{\mathbf{u}}
\newcommand{\bx}{\mathbf{x}}
\newcommand{\be}{\mathbf{e}}
\newcommand{\mH}{\mathcal{H}}

\usepackage{xspace,xcolor}
\usepackage[linecolor=blue!40!white, backgroundcolor=blue!05!white,bordercolor=blue!40!white,textsize=tiny]{todonotes}
\let\origintodo\todo  \newcommand{\xtodo}[2][]{\origintodo[#1]{#2}\xspace}  \let\todo\xtodo

\title{Numerical simulation of the dynamics of molecular markers involved in cell polarisation}

\author{Vincent Calvez\thanks{Unit\'e de Math\'ematiques Pures et Appliqu\'ees, CNRS UMR 5669 \& \'equipe-projet INRIA NUMED, \'Ecole Normale Sup\'erieure de Lyon, 46 all\'ee d'Italie, F-69364 Lyon, France.  ({\tt vincent.calvez@umpa.ens-lyon.fr})} \and Nicolas Meunier\thanks{MAP5, CNRS UMR 8145, Universit\'{e} Paris Descartes, 45 rue des Saints  P\`{e}res
75006 Paris,
France. ({\tt nicolas.meunier@parisdescartes.fr})} \and Nicolas Muller \thanks{MAP5, CNRS UMR 8145, Universit\'{e} Paris Descartes, 45 rue des Saints  P\`{e}res
75006 Paris,
France. {\tt nicolas.muller@parisdescartes.fr})}  \and Raphael Voituriez.
        \thanks{Laboratoire de la mati\`ere condens\'ee, CNRS UMR 7600, Universit\'e Pierre et Marie Curie, 4 Place Jussieu, 75255 Paris Cedex 05 France ({\tt voiturie@lptmc.jussieu.fr})}}

\begin{document}

\maketitle


\begin{abstract}
A cell is polarised when it has developed a main axis of organisation through the reorganisation of its cytosqueleton and its intracellular organelles. Polarisation can occur spontaneously or be triggered by external signals, like gradients of signaling molecules ... Following \cite{Firstpaper} and \cite{Siam_CHMV}, in this work, we study mathematical models for cell polarisation. These models are based on nonlinear convection-diffusion equations. The nonlinearity in the transport term expresses the positive loop between the level of protein concentration localised in a small area of the cell membrane and the number of new proteins that will be convected to the same area. We perform numerical simulations and we illustrate that these models are rich enough to describe the apparition of a polarisome.   
\end{abstract}

%
%
%


\section{Introduction}

Cell polarisation is a major step involved in several important cellular processes such as directional migration, growth, oriented secretion, cell division, mating or morphogenesis. When a cell is not polarised molecular markers (proteins CDC42) are uniformly distributed on the membrane while polarisation is characterized by the concentration of molecular markers in a small area of the cell membrane. In \cite{Altschuler2003}, it has been observed that if the external pheromone concentration is above a critical concentration, polarisation can occur spontaneously. It has also been observed that cell asymmetry can be driven by an external asymmetric stimulus. 

Cell polarisation in yeast cells has been intensively studied during the past decade. Recently, many models describing cell polarisation have been developed. The majority of these models are based on reaction-diffusion systems where polarisation appears as a type of Turing instability \cite{Turing1}, \cite{Turing2}, \cite{Turing3}, or due to stochastic fluctuations \cite{Altschuler}, other models include cytoskeleton proteins as a regulatory factor \cite{Marco2007}, \cite{Altschuler2003}. Many biological studies have shown that the cytoskeleton plays an important role in polarisation. It has been suggested that the cytoskeleton has a positive feedback on molecular markers density. Indeed, disruption of transport along the cytoskeleton greatly reduces the stability of polar cap \cite{Altschuler2003}.  The cell cytoskeleton is a network of long semi-flexible filaments made up of protein subunits \cite{cytoskeleton}. These filaments (mainly actin or microtubules) act as roads along which motor proteins are able to perform a biased ballistic motion and carry various molecules. Molecular markers play a key role in the formation of these filaments.

Following \cite{Firstpaper}, \cite{CRAS} and \cite{Siam_CHMV}, in this work we study models that describe the dynamics of cell polarisation. In these models, molecular markers, such as proteins, diffuse in the cytoplasm and are actively transported along the cytoskeleton. The resulting motion is a biased diffusion regulated by the markers themselves. Using numerical simulations and mathematical heuristics, we observe that the coupling on the velocity field achieves an inhomogeneous distribution of molecular markers without any external asymmetric field. Such an inhomogeneous distribution  is only due to interaction between molecular markers.

Throughout this paper, the density of molecular markers (resp. advection field) is denoted by $\rho(t,\bx)$ (resp. $\bu(t,\bx)$). The advection is obtained though a coupling with the membrane concentration of markers. The cell is figured by the domain $\Omega \subset \mathbb{R}^n$ with $n=1,2$ and a part of the boundary of the domain will be the active membrane denoted by $\Gamma$.  The time evolution of the molecular markers satisfies the following advection-diffusion equation, see \cite{Firstpaper} and \cite{Siam_CHMV}:
\begin{equation}\label{main}
\left\{
\begin{aligned}
\partial_t \rho(t,\bx) & = D\, \Delta \rho(t,\bx) - \chi\, \nabla . \left( \rho(t,\bx) \, \bu (t,\bx) \right), \quad t > 0, \quad \bx \in \Omega, \\
\rho(0,\bx) & = \rho_0 (\bx).
\end{aligned}
\right.
\end{equation}
There is no creation nor degradation of molecular markers in the cell, so the quantity of molecular markers remains constant in time:
\begin{equation}\label{massconservation}
M = \int_{x\in\Omega} \rho_0(\bx) d\bx = \int_{x\in\Omega} \rho(t,\bx) d\bx.
\end{equation}
This condition is ensured by a zero flux boundary condition on the boundary.
A first simplified step is to assume that the cell is essentially bidimensional and to neglect curvature effects. The membrane
boundary is then a 1D line along the $y$-axis and the cytoplasm is parametrized by $\bx=(x,y) \in \mathbb{R}_+ \times \mathbb{R}$. 

The plan of this work is the following. First, we recall the main mathematical results of the simplified model in 1D for $\Omega=(0,\infty)$ and $\Gamma = \{ x = 0 \}$, see \cite{CRAS}, \cite{Siam_CHMV} for more details. Then we study a more realistic model, that includes dynamical exchange of markers on the boundary for a general $\Omega$. This model was introduced in \cite{Firstpaper} and studied in \cite{Siam_CHMV} in the one dimensional case. Here, we will perform a first numerical analysis of this  model in the two dimensional case, for periodic (in one direction) and bounded (in the other direction) domain. Finally, we provide a methodology for parameter estimation by using mathematical heuristics and biological literature.

\subsection{One dimensional case}
In this section, we study the one dimensional case on the half line for $\Omega=(0,\infty)$. The membrane is then the point $\Gamma = \{ x = 0 \}$. For the first model, the advection field towards the membrane is equal to the density of molecular markers on the boundary $\rho(t,0)$. Then we improve this model by considering that only the trapped molecular markers on the membrane contribute to the advection field.
\subsubsection{Simplified model set on the half line}
In \cite{Siam_CHMV} a first mathematical studies  has been done on this model. We define an advection field $\bu(t,x)$ for \eqref{main} $$\bu(t,x)=-\rho(t,0),$$
in such a case \eqref{main} reads as (with $D=1$ and $\chi=1$):
\begin{equation}\label{simplified}
\partial_t \rho(t,x) = \partial_{xx} \rho(t,x) + \rho(t,0) \, \partial_x \rho(t,x), \quad t > 0, \quad x>0,
\end{equation}
with the following zero flux condition on the boundary $\{x=0\}$, that ensures the mass conversation \eqref{massconservation},
\begin{equation}
\partial_x \rho(t,0) + \rho(t,0)^2=0.
\end{equation}
In \cite{Siam_CHMV}, it has been proved that solutions of \eqref{simplified} blow-up in finite time if their masses are above a certain critical mass, $ M> 1 $, and exist globally in time if $ M \leq 1 $. Let us first recall the definition of weak solutions of \eqref{simplified}.
\begin{definition}\label{def:weak}
We say that $\rho(t,x)$ is a weak solution of \eqref{simplified}  on $(0,T)$ if it satisfies:
\begin{equation*}
\rho \in L^\infty(0,T;L^1_+(\R_+))\, , \quad \partial_x \rho \in L^1((0,T)\times \R_+)  \, , \label{eq:flux L1}
\end{equation*}
and $\rho(t,x)$ is a solution of \eqref{simplified} in the sense of distributions in $\mathcal D'(\R_+)$.
\end{definition}

Let us now recall the main results for weak solutions of \eqref{simplified}.
\begin{theorem}[Global existence: $M\leq1$] \label{th:1D} 
Assume that the initial data $\rho_0$ satisfies both $\rho_0 \in L^1(( 1 + x)\D x)$ and $\int_{x>0} \rho_0(x) (\log \rho_0(x))_+ \D x< + \infty$. Assume in addition that $M\leq 1$, then there exists a global weak solution of equation \eqref{simplified}.
\end{theorem}
\begin{theorem}[Blow-up: $M>1$] \label{th:1D BU} Assume $M>1$. Any weak solution of equation \eqref{simplified} with non-increasing initial data $\rho_0$ blows-up in finite time. \\
\end{theorem}
\begin{remark}
It would tempted to interpret blow-up of solutions of the one dimensional model as cell polarisation. But it is to be noticed that concentration of markers on the boundary doesn't mean polarisation. Indeed, consider a radially symmetric 2D cell case. Equation then reduces to the one dimensional one. Above a threshold on the total mass, the convection wins and markers concentrate on the boundary. In some situations, these markers may be homogeneously distributed on the boundary and in such a case there is no symmetry breaking.
\end{remark}

\subsubsection{The model with dynamical exchange of markers at the boundary} 
Such a direct activation of transport on the boundary seems to be unrealistic. Indeed possible occurrence of blow-up in finite time suggests this claim. We improve the previous model by distinguishing between cytoplasmic content $\rho(t,x)$ and the concentration of trapped molecules on the boundary, that will be denoted by $\mu(t)$. The dynamical exchange of markers at the boundary is done with an attachment rate $k_{on}$ and a detachment rate $k_{off}$, hence the time evolution of $\mu(t)$ is
\begin{equation}\label{mu1D}
\dfrac{d}{dt} \mu (t) = k_{on} \, \rho(t,0) - k_{off} \, \mu(t).
\end{equation}
The advection field $\bu(t,x)$ in \eqref{main} is now defined by $$\bu(t,x)=-\mu(t),$$
hence \eqref{main}  (with $D=1$ and $\chi=1$) reads as:
\begin{equation}
\partial_t \rho(t,x) = \partial_{xx} \rho(t,x) + \mu(t) \, \partial_x \rho(t,x), \quad t > 0, \quad x>0,
\end{equation}
with a modified boundary condition
\begin{equation}
\partial_x \rho(t,0) + \rho(t,0)\, \mu(t) = \dfrac{d}{dt} \mu (t).
\end{equation}
This ensures the following mass conservation shared among $\rho(t,x)$ and $\mu(t)$:
\begin{equation}
M = \int_{\mathbb{R}_+} \rho_0(x) dx + \mu_0 = \int_{\mathbb{R}_+} \rho(t,x) dx + \mu(t).
\end{equation}
With equation \eqref{mu1D}, the self-activation of transport by $\rho(t,0)$ is then delayed in time. Since the transport speed is bounded $\mu(t) \leq M$, the solution of the model with dynamical exchange on the boundary exists globally in time. More precisely it is possible, see \cite{Siam_CHMV}, to prove that it converges towards a non trivial stationary state.

\subsection{Two dimensional case : the model with dynamical exchange of markers at the boundary} 
Let $\Omega \subset \mathbb{R}^2$ be the cytoplasm domain, as in the one dimensional case \eqref{mu1D} we consider dynamical exchange of markers at the boundary, so for $\bx \in \Gamma$ we have the evolution in time of $\mu(t,\bx)$
\begin{equation}
\partial_t \mu (t,\bx) = k_{on} \, \rho(t,\bx) - k_{off} \, \mu(t,\bx).
\end{equation}
with a modified boundary condition for $\rho(t,\bx)$ at point $\bx \in \Gamma$
\begin{equation}
(D\, \nabla \rho(t, \bx) - \chi \, \rho(t,\bx)\, \bu(t,\bx) ). \vec n_\bx  = - \partial_t \mu (t,\bx).
\end{equation}
where $\vec n_\bx$ is the outward normal to $\Gamma$.
This ensures the following mass conservation sharing by $\rho(t,\bx)$ and $\mu(t,\bx)$:
\begin{equation}
M = \int_{\Omega} \rho_0(\bx) d\bx + \int_{\Gamma} \mu_0(\bx) d\bx = \int_{\Omega} \rho(t,\bx) d\bx + \int_{\Gamma} \mu(t,\bx) d\bx.
\end{equation}
We consider the following advection field deriving from a harmonic potential modeling the transport by actin filaments (cytoskeleton): 
\begin{equation}\label{c}
\bu(t,\bx) = \nabla c(t,\bx), \mbox{ where } \begin{cases}- \Delta c(t,\bx) = 0, & \mbox{ if } \bx \in \Omega, \\
\nabla c(t,\bx) . \vec n_\bx = S(\bx) \mu(t,\bx), & \mbox{ if } \bx \in \Gamma.
\end{cases}
\end{equation}
This advection field orientation is due to the actin networks.
\begin{center}
\includegraphics[scale=0.7]{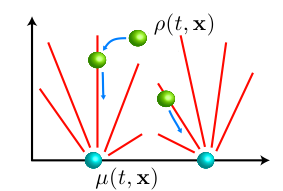}
\end{center}
Actin filaments are attached on the membrane and randomly distributed, there orientations are mixed up. We also add the external pheromone concentration at $\bx \in \Gamma$ which acts by the mating-pheromone MAPK cascade on the actin transport.

In dimension 2, we have global existence for the model without exchange on the boundary (replacing equation \eqref{c} by $\nabla c(t,\bx) . \vec n_\bx = S(\bx) \rho(t,\bx) \mbox{ if } \bx \in \Gamma$) with $\Omega = (0,+\infty) \times \mathbb{R}$ and $\Gamma = \{0\} \times \mathbb{R}$. For clarity, we recall this result,  see \cite{Siam_CHMV} for more details.
\begin{theorem}[Global existence in dimension 2]\label{thdim2}
Assume that the advection field satisfies the two following conditions: $\nabla\cdot \bu \geq 0$ and $\bu(t,0,y)\cdot \be_e = \rho(t,0,y)$. Assume that the initial data $\rho_0$ satisfies both $\rho_0 \in L^1(( 1 + |\bx|^2)\D \bx)$ and $\|\rho_0\|_{L^2}$ is smaller than some constant $c$. Then there exists a global weak solution to equations \eqref{main}-\eqref{massconservation}. \\ 
\end{theorem}

In the two dimensional case, for the model with exchange on the boundary,  blow-up or global existence have not been proved yet. In this work, we make a first step in this direction by using a mathematical heuristic and numerical simulations.

\subsection{Heuristics}
The mathematical analysis performed in \cite{Siam_CHMV} has demonstrated that a class of models exhibit pattern formation (either blow-up or convergence towards a non homogeneous steady state) under some conditions. However the main question still remains unanswered: do these models describe cell polarisation or not? Thus in order to provide a first answer to this question, we will perform numerical simulations. Our aim is to see if, under some conditions, the model leads to a concentration of markers, not only on the boundary, but on a small region of the boundary. In such a case polarisation occurs. In order to obtain more information on the critical value distinguishing the polarised case and the stable case, in the two dimensional case we will use a mathematical heuristics that we describe now.  
Let $\bx = (x,y) $ be in $\Omega = \mathbb{R}_+ \times \mathbb{R}$ and let  $\Gamma=\{0\} \times \mathbb{R}$ be the boundary, we have that
\begin{equation}
\bu(t,\bx) = \nabla c(t,\bx), \mbox{ where } \begin{cases}-\Delta c(t,\bx) = 0, & \mbox{ if } \bx \in \mathbb{R}_+ \times \mathbb{R}, \\
-\partial_x c(t,0,y)= S(y) \mu(t,y), & \mbox{ if } y \in \mathbb{R},
\end{cases}
\end{equation}
hence,  see  \cite{Evans} e.g., it is well known that 
$$c(x,y) = - \frac{1}{\pi} \int_{y'\in \mathbb{R}}\log (\sqrt{(y-y')^2+x^2}) (S \mu)(y') dy'.$$
The tangential component at the boundary is then given by
\begin{equation*}
\bu (t,0,y) \cdot \vec \be_y = -\mH ( S \mu) (y) \, , \quad y \in \mathbb{R},
\end{equation*}
where $\mH$ denotes the one-dimensional Hilbert transform that we recall now, see  \cite{KellerSegelHilbert} e.g., with respect to the $y$ variable:
\begin{equation*}
\mH (f\mu) (y) = \dfrac1\pi{\rm p.v.}\int_{\R}\dfrac1{y-y'} f(y')\, dy'.
\end{equation*}
Integrating the main equation \eqref{main} with respect to $x$ with zero flux condition on $\Gamma = \{x=0\}$, we obtain:
\begin{equation*}
\partial_t \int_{x>0} \rho(t,x,y)\, dx  = D \, \partial_{yy} \left(\int_{x>0} \rho(t,x,y)\, dx\right) - \chi \, \partial_y \left(\int_{x>0} \rho(t,x,y) (\bu(t,x,y) \cdot \vec \be_y) dx \right).
\end{equation*}
In the super-critical case, numerical simulations, see \cite{NicoPhD}, suggest that the density $\rho(t,\bx)$ concentrates on the boundary $\{x = 0\}$. Assuming $\rho(t,x,y) = \nu(t,y)\delta(x = 0)$, we can formally write the dynamics of $\nu(t,y)$ as follows:
\begin{equation*}
\partial_t \nu(t,y)  = D \, \partial_{yy} \nu(t,y) + \chi  \, \partial_y \left(\nu(t,y) \mH ( S \mu) (y) \right).
\end{equation*}
Assuming $S$ constant on $\mathbb{R}$ and $\mu(t,y) = \frac{k_{on}}{k_{off}} \nu(t,y)$ for $y \in \mathbb{R}$, it reads as
\begin{equation*}
\partial_t \nu(t,y)  = D \, \partial_{yy} \nu(t,y) + \chi   S \frac{k_{on}}{k_{off}} \, \partial_y \left(\nu(t,y) \mH ( \nu) (y) \right)\, .
\end{equation*}

Hilbert transform has a critical singularity to offset the diffusion on this equation \cite{KellerSegelHilbert}. We have a blow-up if $\int_\R \nu(t,y)\, dy = M$ is above $\frac{2 \pi D k_{off}}{ S \chi k_{on}}$. This is a first step to observe a critical mass phenomenon and this may lead to blow-up if the mass is large enough. In this way, we define an order of magnitude for some parameters. 

It is to be noticed that this latter criterion is valid for an infinite domain, namely $y\in \R$. In the case of a cell, the domain will be finite and the existence of such a dichotomy has not been proved yet. In order to see if such a dichotomy holds true we will perform numerical simulations. This is the object of the following section.

\section{Numerical analysis}

We first give a discretization of the convection-diffusion model set on a 1D periodic domain. This first step allows us introducing the discretization of this model on a 2D domain which is periodic in one direction and bounded on the other direction.

\subsection{One dimensional case} Let $u(t,x)$ be a given function. We consider the following advection-diffusion equation on the periodic domain $\Omega = \mathbb{R}/\mathbb{Z}$
\begin{equation}\label{1D}
\partial_t \rho(t,x) =  \partial_{x}  (\partial_{x} \rho(t,x) - u(t,x) \, \rho(t,x)), \quad t > 0, \quad x\in \Omega.
\end{equation}
Let $t^n = n\, dt$ be the time discretization and $\{x_j=j \, dx, j \in \{1,...,N_x\}\}$ be the space discretization of the periodic interval $\mathbb{R}/\mathbb{Z}$. Since the equations of the model are written in a conservative form, the natural framework to be used for the spatial discretization is the finite volume framework.  We hence introduce the control volume defined for $j \in \{1,...,N_x\}$
\begin{equation}
V_j = (x_{j-\frac{1}{2}},x_{j+\frac{1}{2}}).
\end{equation}
Let $\rho_j^n$ (resp. $u^n_{j+\frac{1}{2}}$) be the approximated value of the exact solution $\rho(t^n,x_j)$ (resp. $u(t^n,x_{j+\frac{1}{2}})$), the classical upwind scheme for \eqref{1D} reads as
\begin{eqnarray}
\frac{\rho_j^{n+1} - \rho_j^n}{dt} = \frac{\mathcal{F}_{j+\frac{1}{2}} - \mathcal{F}_{j-\frac{1}{2}}}{dx}, \quad j \in \{1,...,N_x\},
\end{eqnarray}
where the numerical flux $\mathcal{F}_{j+\frac{1}{2}}$ and $\mathcal{F}_{j-\frac{1}{2}}$ are defined by
\begin{eqnarray*}
\mathcal{F}_{j+\frac{1}{2}} = \frac{\rho_{j+1}^{n+1} -  \rho_j^{n+1}}{dx}- \textit{A}^{up}(u_{j+\frac{1}{2}}^n,\rho_j^n,\rho_{j+1}^n), \\
\mathcal{F}_{j-\frac{1}{2}} = \frac{\rho_{j}^{n+1} -  \rho_{j-1}^{n+1}}{dx}- \textit{A}^{up}(u_{j-\frac{1}{2}}^n,\rho_{j-1}^n,\rho_j^n),
\end{eqnarray*}
with the advection numerical flux is given by
\begin{equation}\label{Aup}
\textit{A}^{up}(u,x_-,x_+) = \begin{cases}
u \, x_-, \quad \mbox{ si } u > 0, \\ 
u \, x_+, \quad \mbox{ si } u < 0.
\end{cases}
\end{equation}
The periodic flux condition on boundary reads as $\mathcal{F}_{ \frac{1}{2}} = \mathcal{F}_{N_x + \frac{1}{2}}$ and we set the value $u^n_{ \frac{1}{2}} = u^n_{N_x + \frac{1}{2}}$.
The diffusion part is treated implicitly and it is then unconditionally stable, while the advection term is treated explicitly. The CFL condition of the scheme is $$\left|\left| \left(u^n_{ j+\frac{1}{2}}\right)_{j \in \{1,...,N_x\}} \right|\right|_\infty < \frac{dx}{dt}.$$
We define the column vector $\rho^{n} = \begin{pmatrix} \rho_1^n & \rho_2^n & \dots & \rho_{N_x}^n  \end{pmatrix}^T$. As usual, see e.g. \cite{Allaire}, the discrete heat matrix $A \in M_{N_x} (\mathbb{R})$ with periodic flux condition on the boundary is defined as 
\begin{equation}\label{defA}
A = 
\begin{pmatrix} 
2 + \frac{dx^2}{dt} & -1 & & & -1 \\ -1 & 2 + \frac{dx^2}{dt} & \ddots \\ & \ddots & \ddots & \ddots \\& & \ddots & 2 + \frac{dx^2}{dt} & -1 \\ -1 & & & -1 & 2 + \frac{dx^2}{dt} 
\end{pmatrix}.
\end{equation}
Periodic flux condition adds the top right term and the bottom left term. Next, in order to use $A^{up}$ defined by equation \eqref{Aup}, we define $(u)^+=\max(u,0)$ and $(u)^-=\min(u,0)$. The discrete advection matrix $B \in M_{N_x} (\mathbb{R})$ with periodic flux condition on the boundary is then defined as in \cite{Allaire}
\begin{equation}\label{defB}
B =  \frac{dx^2}{dt} I_{N_x} - dx  \,
\begin{pmatrix} 
\left(u_{\frac{3}{2}}^n\right)^+ & \left(u_{\frac{3}{2}}^n\right)^- & & &   \\ & \ddots & \ddots  \\ & & \left(u_{j+\frac{1}{2}}^n\right)^+ & \left((u_{j+\frac{1}{2}}^n\right)^-  \\ & & & \ddots & \left(u_{N_x-\frac{1}{2}}^n\right)^- \\ \left(u_{N_x+\frac{1}{2}}^n\right)^- & & & & \left(u_{N_x+\frac{1}{2}}^n\right)^+ 
\end{pmatrix}
\end{equation}
\begin{equation*}
+ dx  \,
\begin{pmatrix} 
\left(u_{\frac{1}{2}}^n\right)^- & & & & \left(u_{\frac{1}{2}}^n\right)^+  \\ \left(u_{\frac{3}{2}}^n\right)^+ & \ddots &   \\ & \left(u_{j-\frac{1}{2}}^n\right)^+ & \left(u_{j-\frac{1}{2}}^n\right)^- &  \\ & & \ddots & \ddots &  \\  & & & \left(u_{N_x-\frac{1}{2}}^n\right)^+ & \left(u_{N_x-\frac{1}{2}}^n\right)^-
\end{pmatrix}.
\end{equation*}
We use a standard numerical method to invert the symmetric positive definite matrix $A$. Finally,  at each time step we resolve
\begin{equation*}
\rho^{n+1}= A^{-1} \, B \, \rho^{n}.
\end{equation*}

\subsection{Two dimensional case}
We  perform numerical simulations on the model with dynamical exchange of markers at the boundary. In this work, we assume that the cell occupies a disk of radius $r>0$. Furthermore for simplicity, we consider a bounded-periodic domain $\Omega = [0,r] \times \mathbb{R}/2\pi r\mathbb{Z}$ with $\Gamma = \{r\} \times \mathbb{R}/2\pi r\mathbb{Z}$. This simplifies our numerical approach by using finite difference schemes on Cartesian grid. We start with the numerical study of the equation on $\rho$ by assuming that the advection field $u(t,\bx) = \nabla c(t,\bx)$ is known. Then we perform the discretization of $c$.

In this section, for simplicity we fix all the parameters values to 1 except $M$. Let us first recall the model with dynamical exchange of markers at the boundary on $\Omega = [0,r] \times \mathbb{R}/2\pi r\mathbb{Z}$:
\begin{eqnarray}
\partial_t \rho = \nabla . \left( \nabla \rho -  \rho \, \nabla c \right), \mbox{ in } (0,r) \times \mathbb{R}/2\pi r\mathbb{Z}, \label{rho1}\\
\partial_x \rho -  \rho \, \partial_x c =  -\partial_t \mu, \mbox{ on } \{r\} \times \mathbb{R}/2\pi r\mathbb{Z}, \label{rho2} \\
\partial_x \rho -  \rho \, \partial_x c =  0, \mbox{ on } \{0\} \times \mathbb{R}/2\pi r\mathbb{Z} \label{rho3}.
\end{eqnarray}
Dynamical exchange markers on active boundary $\{r\} \times \mathbb{R}/2\pi r\mathbb{Z}$ is given by
\begin{eqnarray}
\partial_t \mu=\rho-\mu, \mbox{ on } \{r\} \times \mathbb{R}/2\pi r\mathbb{Z}\label{mu2D}.
\end{eqnarray}
Laplace equation on $c$ with non appropriate Neumann conditions on a bounded domain is ill-posed,  see \cite{Allaire} e.g. In order to handle this problem, we add a degradation term:
\begin{eqnarray}
-\Delta c + \alpha \, c = 0, \mbox{ in } (0,r) \times \mathbb{R}/2\pi r\mathbb{Z}, \label{c1}\\
- \partial_x c =  \, \mu, \mbox{ on } \{r\} \times \mathbb{R}/2\pi r\mathbb{Z}, \label{c2} \\
- \partial_x c = 0, \mbox{ on } \{0\} \times \mathbb{R}/2\pi r\mathbb{Z} \label{c3}.
\end{eqnarray}
We take random initial conditions $c$, $\mu_0$ and $\rho_0$ satisfying the following mass conservation
\begin{eqnarray}\label{massnum}
\int_\Omega \rho_0 + \int_\Gamma \mu_0 = M.
\end{eqnarray}
Let $t^n = n\, dt$ be the time discretization and $\{x_j=j \, dx, j \in \{1,...,N_x\}\}$ (resp. $\{y_k=k \, dy, k \in \{1,...,N_y\}\}$) be the space discretization of the bounded interval $[0,r)$ (resp. periodic interval $\mathbb{R}/ 2 \pi r \mathbb{Z}$). We introduce the control volume $W_{(j,k)}  \subset \mathbb{R}^2$
\begin{equation*}
W_{(j,k)} = (x_{j-\frac{1}{2}},x_{j+\frac{1}{2}}) \times (y_{k-\frac{1}{2}},y_{k+\frac{1}{2}}).
\end{equation*}

Let $P^n_{(j,k)}$ (resp. $\mu^n_{k}$) be the approximated value of the exact solution $\rho(t^n,x_j,y_k)$ of equations \eqref{rho1},\eqref{rho2},\eqref{rho3} and \eqref{massnum} (resp. $\mu(t^n,y_k)$ of equations \eqref{mu2D} and \eqref{massnum}). 
Let  $c_{(j,k)}$ be the approximated value of the exact solution $c(x_j,y_k)$ of equations \eqref{c1},\eqref{c2} and \eqref{c3}. 

\subsubsection{Equation on $\mu$}
We can resolve at each time step for $k \in \{1,...,N_y\}$
\begin{eqnarray*}
\mu^{n+1}_k=\mu^{n}_k + dt \, (\rho^n_k-\mu^n_k).
\end{eqnarray*}

\subsubsection{Equation on $c$}
For simplicity, we call $\mathcal{F}$ the numerical flux as in the 1D case, we can write the following scheme: for $(j,k) \in \{1,...,N_x\} \times \{1,...,N_y\}$
\begin{eqnarray*}
\frac{\mathcal{F}_{(j+\frac{1}{2},k)} - \mathcal{F}_{(j-\frac{1}{2},k)}}{dx} + \frac{\mathcal{F}_{(j,k+\frac{1}{2})} - \mathcal{F}_{(j,k-\frac{1}{2})}}{dy}  - \alpha c_{(j,k)} = 0.
\end{eqnarray*}
with numerical flux defined by
\begin{eqnarray*}
\mathcal{F}_{(j+\frac{1}{2},k)} = \frac{c_{(j+1,k)} - c_{(j,k)} }{dx},\\
\mathcal{F}_{(j-\frac{1}{2},k)} = \frac{c_{(j,k)} - c_{(j-1,k)} }{dx}.
\end{eqnarray*}
The zero flux boundary condition \eqref{c3} impose that $\mathcal{F}_{(\frac{1}{2},k)} = 0$ and the boundary condition \eqref{c2} $\mathcal{F}_{(N_x + \frac{1}{2},k)} = - \mu^n_k$ for $k \in \{1,...,N_y\}$. Similarly, the periodic conditions impose $\mathcal{F}_{(j,N_y+\frac{1}{2})} = \mathcal{F}_{(j,\frac{1}{2})}$  for $j \in \{1,...,N_x\}$.
We define the column vector $\mathcal{C}$ by $\mathcal{C}(k+(j-1) N_y) = C_{(j,k)}$ with $(j,k) \in \{1,...,N_x\} \times \{1,...,N_y\}$.
As previously the rigidity matrix $A_{2D,\alpha}$ is defined by
\begin{equation*}
A_{2D,\alpha} =
\begin{pmatrix} 
A_{\alpha}+Id & -Id &  & & \\ -Id & A_{\alpha}+2 \, Id & -Id &  & & \\ & -Id & A_{\alpha}+2 \, Id & \ddots & &  \\ & & \ddots & \ddots & \ddots & & \\  & & & \ddots & A_{\alpha}+2 \, Id & -Id & \\ & & & & -Id & A_{\alpha}+2 \, Id & -Id \\ & & & & & -Id & A_{\alpha}+Id
\end{pmatrix},
\end{equation*}
where discrete Poisson matrix $A_{\alpha} \in M_{N_y} (\mathbb{R})$ in 1D with periodic flux condition on boundary is defined by
\begin{equation*}
A_{\alpha} = 
\begin{pmatrix} 
2 + \alpha \, dx^2 & -1 &  & & -1 \\ -1 & 2 + \alpha \, dx^2 & \ddots &  &  \\ & \ddots & \ddots & \ddots & &  \\ & & \ddots  & 2 + \alpha \, dx^2 & -1 &  \\ -1 & & & -1 & 2 + \alpha \, dx^2
\end{pmatrix}.
\end{equation*}
The flux boundary condition $\{r\} \times \mathbb{R}/2\pi r\mathbb{Z}$ imposes this right hand side column vector of length $N_x \, N_y$:
\begin{equation*}
R_c = - dx \begin{pmatrix} (\mu^n_k)_{k} & 0 & \dots & 0 \end{pmatrix}.
\end{equation*}
We use a standard numerical method to invert the symmetric positive definite matrix $A_{2D,\alpha}$ and then resolve at each time step
\begin{equation*}
\mathcal{C} = A_{2D,\alpha}^{-1} \, R_c.
\end{equation*}

\subsubsection{Equation on $\rho$}
For simplicity, we call $\mathcal{F}$ the numerical flux as in the previous cases, we can write the upwind scheme as follows:
\begin{eqnarray*}
\frac{P_{(j,k)}^{n+1} - P_{(j,k)}^n}{dt} & = & \frac{\mathcal{F}_{(j+\frac{1}{2},k)} - \mathcal{F}_{(j-\frac{1}{2},k)}}{dx}  + \frac{\mathcal{F}_{(j,k+\frac{1}{2})} - \mathcal{F}_{(j,k-\frac{1}{2})}}{dy},
\end{eqnarray*}
where $u^n_{(j+\frac{1}{2},k)} = \frac{c_{(j+1,k)}^{n}-c_{(j,k)}^{n}}{dx}$, $u^n_{(j-\frac{1}{2},k)} = \frac{c_{(j,k)}^{n}-c_{(j-1,k)}^{n}}{dx}$ and the numerical flux is defined by
\begin{eqnarray*}
\mathcal{F}_{(j+\frac{1}{2},k)} = \frac{\rho_{(j+1,k)}^{n+1} - \rho_{(j,k)}^{n+1} }{dx} - \textit{A}^{up}\left(u^n_{(j+\frac{1}{2},k)},P_{(j,k)}^{n}, P_{(j+1,k)}^{n} \right), \\
\mathcal{F}_{(j-\frac{1}{2},k)} = \frac{\rho_{(j,k)}^{n+1} - \rho_{(j-1,k)}^{n+1} }{dx} - \textit{A}^{up}\left(u^n_{(j-\frac{1}{2},k)},P_{(j-1,k)}^{n}, P_{(j,k)}^{n} \right).
\end{eqnarray*}
The zero flux boundary conditions \eqref{rho3} impose $\mathcal{F}_{(\frac{1}{2},k)} = 0$, while the flux boundary conditions \eqref{rho2} impose $\mathcal{F}_{(N_x + \frac{1}{2},k)} = -\frac{\mu_k^{n+1}-\mu_k^n}{dt}$ for $k \in \{1,...,N_y\}$. 
Similarly, the periodic conditions impose $\mathcal{F}_{(j,N_y+\frac{1}{2})} = \mathcal{F}_{(j,\frac{1}{2})}$  for $j \in \{1,...,N_x\}$.
We define the column vector $\mathcal{P}^n$ by $\mathcal{P}^n (k+(j-1) N_y) = P^n_{(j,k)}$ with $(j,k)\in \{1,...,N_x\} \times \{1,...,N_y\}$. In what follows for simplicity we consider that $dx=dy$. We define the rigidity matrix $A_{2D} \in M_{N_x N_y} (\mathbb{R})$ with $A\in M_{N_y} (\mathbb{R})$ defined by \eqref{defA}:
\begin{equation*}
A_{2D} =
\begin{pmatrix} 
A + Id& -Id &  & & \\ -Id & A + 2 \, Id & \ddots &  & & \\ &  \ddots & \ddots & \ddots & & \\  &  & \ddots & A + 2 \, Id  & -Id & \\ & & & -Id & A + Id 
\end{pmatrix}.
\end{equation*}
We define the following diagonal matrices for $j \in \{1,...,N_x\}$,  $U^{+}_{j+\frac{1}{2}} \in M_{N_y}(\mathbb{R})$ and $U^{-}_{j+\frac{1}{2}} \in M_{N_y}(\mathbb{R})$:
\begin{eqnarray*}
U^{+}_{j+\frac{1}{2}} =
\begin{pmatrix} 
\ddots & & & &  \\ & (u_{(j+\frac{1}{2},k-1)}^n)^+ &   \\ & & (u_{(j+\frac{1}{2},k)}^n)^+  &  &  \\ & & & (u_{(j+\frac{1}{2},k+1)}^n)^+ &   \\  & & & & \ddots
\end{pmatrix}, \\
U^{-}_{j+\frac{1}{2}} =
\begin{pmatrix} 
\ddots & & & &  \\ & (u_{(j+\frac{1}{2},k-1)}^n)^- &   \\ & & (u_{(j+\frac{1}{2},k)}^n)^-  &  &  \\ & & & (u_{(j+\frac{1}{2},k+1)}^n)^- &   \\  & & & & \ddots
\end{pmatrix}.
\end{eqnarray*}
With $B\in M_{N_y} (\mathbb{R})$ defined by equation \eqref{defB}, the discrete advection matrix $B_{2D}\in M_{N_x N_y} (\mathbb{R})$ with zero flux boundary condition in the $x$-axis direction and periodic flux boundary condition in the $y$-axis direction is defined by 
\begin{eqnarray*}
B_{2D} & = & \begin{pmatrix} 
B& &  & & \\ & B &  &  & & \\ &   & \ddots &  & & \\  &  &  & B &  & \\ & & & & B
\end{pmatrix}
- dx  \,
\begin{pmatrix} 
U^{+}_{\frac{3}{2}} & U^{-}_{\frac{3}{2}}  \\ & \ddots & \ddots  \\ & & U^{+}_{j+\frac{1}{2}} & U^{-}_{j+\frac{1}{2}}  \\ & & & \ddots & U^{-}_{N_x-\frac{1}{2}} \\  & & & & U^{+}_{N_x+\frac{1}{2}}
\end{pmatrix} \\
& + & dx  \,
\begin{pmatrix} 
U^{-}_{\frac{1}{2}} & & & &  \\ U^{+}_{\frac{3}{2}} & \ddots &   \\ & U^{+}_{j-\frac{1}{2}}  & U^{-}_{j-\frac{1}{2}} &  \\ & & \ddots & \ddots & \\ & & & U^{+}_{N_x-\frac{1}{2}}  & U^{-}_{N_x-\frac{1}{2}}
\end{pmatrix}.
\end{eqnarray*}
The flux boundary condition $\{r\} \times \mathbb{R}/2\pi r\mathbb{Z}$ imposes this right hand side column vector of length $N_x \, N_y$:
\begin{equation*}
R_\rho = - dx \begin{pmatrix} (\frac{\mu^{n+1}_k-\mu^n_k}{dt})_{k} & 0 & \dots & 0 \end{pmatrix}.
\end{equation*}
We use a standard numerical method to inverse the symmetric positive definite matrix $A_{2D}$ and then resolve at each time step
\begin{equation*}
\mathcal{P}^{n+1} = A_{2D}^{-1} \, (B_{2D} \mathcal{P}^n+R_\rho).
\end{equation*}

\subsection{Graphics}

With the previous numerical analysis, we implement all numerical simulations using MATLAB. We test different values of $M$:
\begin{center}
\begin{figure}[!h]
\begin{tabular}{ll}
\textbf{A)}   &  \textbf{B)}  \\
 \includegraphics[scale=.3]{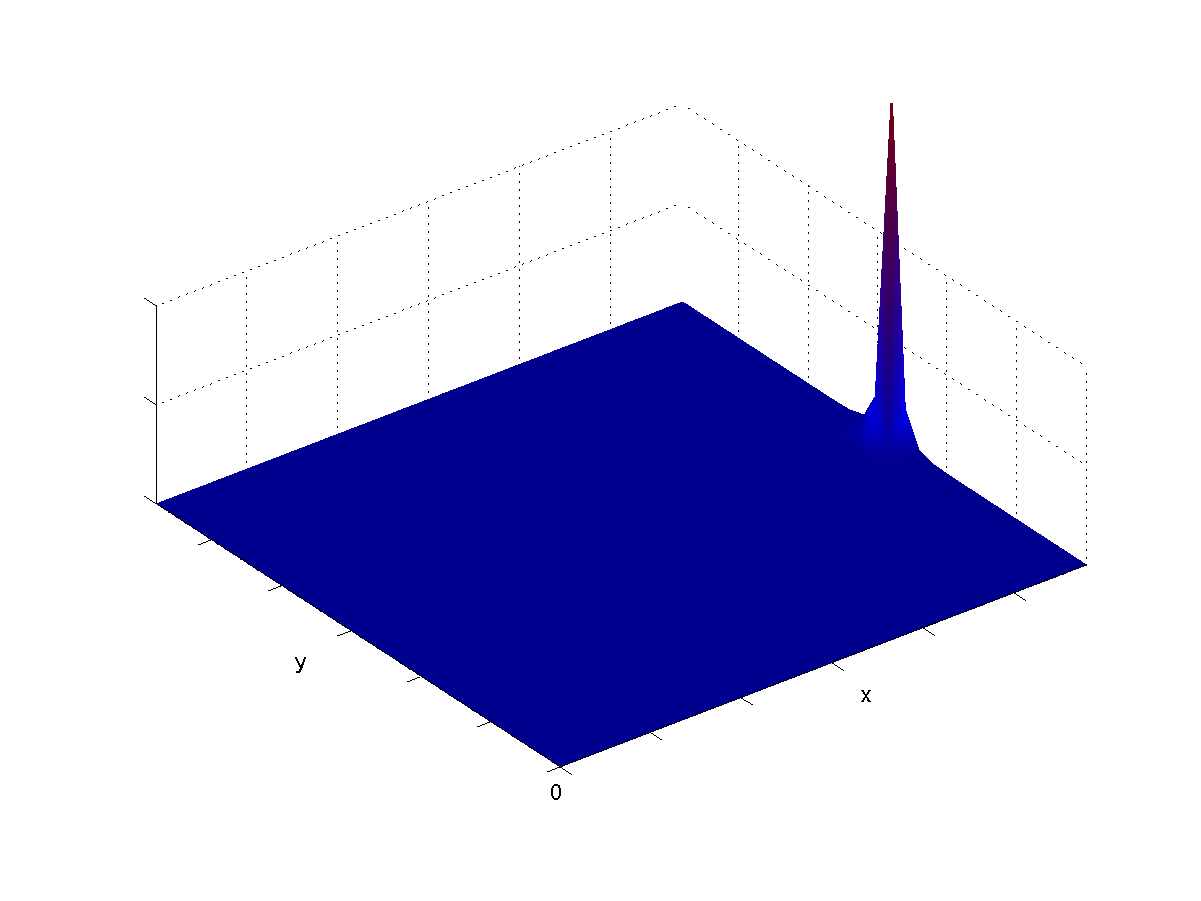}   & \includegraphics[scale=.3]{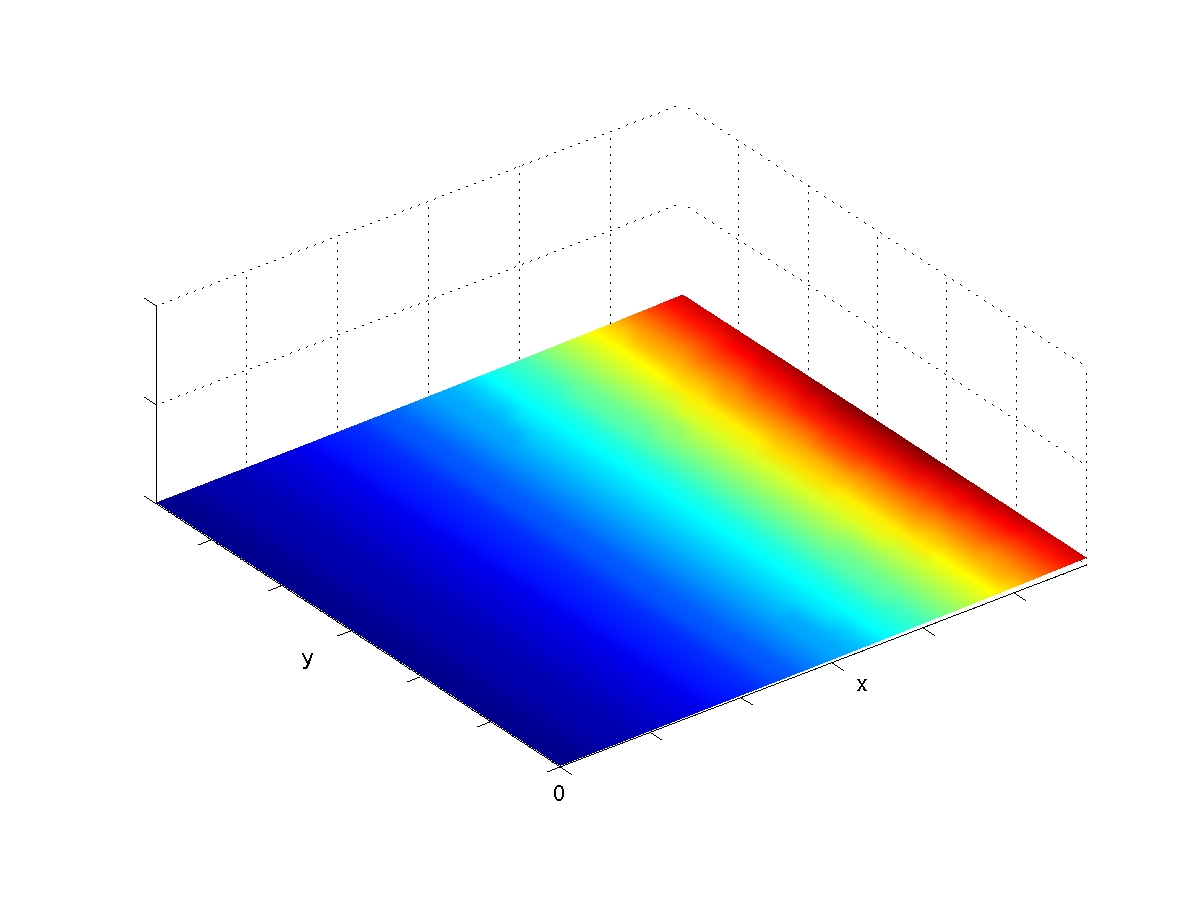}
\end{tabular}
 \caption{Numerical Simulations on $\Omega = [0,1] \times \mathbb{R}/2\pi \mathbb{Z}$ with $\Gamma = \{r\} \times \mathbb{R}/2\pi r \mathbb{Z}$ and all parameters equal to 1. A) For $M=20$ greater enough, symmetry breaking appears. Molecular markers are concentrated on one point of the membrane in finite time. B) For $M=0.01$ small, steady state is homogeneous in the $y$-axis.}\label{fig1}
\end{figure}
\end{center}


\section{Conclusion}

In this work we have provided a first answer to the following question: do the nonlinear convection-diffusion models given in \cite{Firstpaper}  and \cite{Siam_CHMV} describe cell polarisation or not?  To do so we have used both a mathematical heuristic and numerical simulations. Numerical simulations were necessary because the heuristic is only valid for an infinite geometry while the cell is obviously finite.  The numerical simulations ensure that solutions develop symmetry breaking over a critical value $M^*$ given us a first justification of the mathematical heuristic. In a further work, we will estimate an approximate value of this critical mass.
 

\medskip

\bibliographystyle{Siam}

\end{document}